\documentclass[11pt, oneside]{article}   	
\usepackage{geometry}                		
\usepackage{graphicx}				
\usepackage{amssymb}
\usepackage{amsmath}
\usepackage{amsthm}
\usepackage{wasysym}

\usepackage{xcolor}

\usepackage[pdftex]{hyperref}

\usepackage{cleveref}
\usepackage{autonum}

\numberwithin{equation}{section}

\theoremstyle{plain}
\newtheorem{theorem}{Theorem}

\theoremstyle{definition}
\newtheorem{definition}[theorem]{Definition}

\theoremstyle{remark}
\newtheorem{remark}[equation]{Remark}
\newtheorem{example}[theorem]{Example}

\newcommand{\oeis}[1]{\color{blue}\textnormal{\href{https://oeis.org/#1}{#1}}\color{black}\ in \cite{Sloane16}}

\title{Symbolic dynamical scales: \\ modes, orbitals, and transversals}

\author{Ricardo G\'omez A\'iza}

\begin{document}

\maketitle

\begin{abstract}
We study classes of musical scales obtained from shift spaces in symbolic dynamics through the first symbol rule, which yields scales in any $n$-TET tuning system. The modes are thought as elements of orbit equivalence classes of cyclic shift actions on languages, and we study their orbitals and transversals. We present explicit formulations of the generating functions that allow us to deduce the orbital and transversal dimensions of classes of musical scales generated by vertex shifts, for all $n$, in particular for the 12-TET tuning system.
\end{abstract}





\section{Introduction and main result}


A \emph{symbolic sequential scale} is a musical scale obtained
from a (mathematical) symbolic sequence, according to certain (mathematical) rule.
For example, the \emph{standard Thue-Morse scales}, introduced in \cite{GomezNasser20},
are defined as certain sets of scales obtained from the
Thue-Morse binary sequence, with coding rule the 
binary representation of scales.
This rule can be generalized to sequences over larger alphabets,
for example as the \emph{first symbol rule}, that likewise
generates scales on every $n$-TET tuning system;
it consists on coupling a block that occurs in the sequence with a
``symbolic chromatic scale'' of the same length, that is,
a block formed with increasingly ordered elements of a set of ``symbolic notes'',
and then letting the scale
defined by the rule be formed with the notes carrying the distinguished
initial symbol of the corresponding block. 
The rule can be applied to sets of symbolic sequences.
Our goal here is to present a more general formalism to study
symbolic sequential scales of this type, with the
number of notes on each scale as a parameter, 
together with their modes,
using techniques from both symbolic dynamics \cite{LindMarcus95}
and analytic combinatorics \cite{FlajoletSedgewick09}.

We consider \emph{shift spaces}
as sets of symbolic sequences, that is, 
sets of (bi)infinite sequences that avoid 
a given set of (finite) \emph{forbidden} configurations,
and they come with a natural $\mathbb Z$-action $\sigma$ by (left) translation.
A shift space is characterized by its \emph{language},
which is defined as the union of the \emph{admissible} configurations
that occur in its sequences. Thus, 
a shift space $X$ can be constructed
by specifying either its language $\mathcal L(X)$ (e.g., as the admisible
blocks in a sequence like the Thue-Morse sequence),
or a set of forbidden configurations
(e.g., the distinguished symbol rule will certainly yield admissible scales with
no half tones between consecutive notes if the forbidden set
contains the cartesian product of the alphabet).
\emph{Shifts of finite type} (SFTs) are shift spaces defined
by \emph{finite} forbidden sets, and they always conjugate
to \emph{vertex shifts}, which are shift spaces
defined by forbidden sets formed by blocks of length two.
Vertex shifts are well understood symbolic dynamical systems,
they posses matrix representations that provide algebraic and
analytic tools to study dynamic properties, like entropy, periodic
points and their zeta functions, etc.

The modes of symbolic sequential scales are 
instances of orbit equivalence classes of a
\emph{cyclic} shift action $\alpha$ on \emph{finite} sequences
over some (countable) alphabet $\mathcal S$.
\emph{Orbitals} are unions of $\alpha$-orbit equivalence classes,
and they are the subsets upon which $\alpha$ acts.
In general, arbitrary subsets of the full language $\mathcal S^*$
are not orbitals, like the language of a
shift space
(in fact, the language of a shift space is an orbital if and only if
the space is
a full shift).
Likewise,
if the distinguished symbol rule is formalized
as a block function $\varphi \colon \mathcal L(X) \to \mathcal S^*$
valued on finite symbolic sequences over some alphabet $\mathcal S$
that represent musical scales
in a way that the $\alpha$-orbits correspond to the modes of the scales
(for eample, see \eqref{eq:InducedScale}),
then $\varphi\big(\mathcal L (X)\big)$ is not, in general, an orbital.
Thus we consider orbitals
generated by subsets
$B \subseteq \mathcal S^*$ as unions of the $\alpha$-orbit equivalence
classes of their elements, and a \emph{transversal} is a set of
representatives of the generated orbital.

\begin{remark}
\label{remark:problem}
If $B$ is a set of musical scales, like $\varphi\big(\mathcal L (X)\big)$,
then the cardinality of a transversal
is the number of ``essentially different'' scales an instrumentalist would have to learn to play
any scale in $B$, together with all its modes, for a total number of scales that corresponds to the
cardinality of its generated orbital.
\end{remark}

We will refer to these cardinalities as \emph{transversal}
and \emph{orbital dimensions}. Since any set $B$ decomposes into
a sequence $(B_n)_{n\geq 0}$ with $B_n\triangleq B\cap \mathcal S^n$,
there are \emph{transversal} and \emph{orbital}
generating functions $\textnormal{dim}_{\mathsf T}^{B}(z)$
and $\textnormal{dim}_{\mathsf O}^{B}(z)$, respectively.
Thus we aim to find 
transversal and orbital generating functions of classes of
musical scales generated by shift spaces.
We use integer compositions as the combinatorial model for
the class of all musical scales,
in particular because their $\alpha$-orbits represent the modes of the scales (as \emph{wheels}).
Integer compositions are represented as sequences of positive integers,
their generating functions, including bivariate versions marking several
parameters like the number of summands, are well understood \cite{FlajoletSedgewick09}.
With them (see Theorem \ref{thm:ScalesCompositions}),
and the interplay of the $\sigma$-action on sequences and
the $\alpha$-action on languages, 
it is possible to make the main formulation
that is required to deduce
transversal and orbital generating functions
of musical scales generated by 
vertex shifts (with reference to Remark \ref{remark:problem}):

\begin{theorem}
\label{thm:main}
	Let $X \subseteq \mathcal A^{\mathbb Z}$ be an irreducible vertex shift and choose a
	symbol $\mathfrak s\in\mathcal A$.
	Then there is a set $\mathcal K(\mathfrak s)\subseteq \mathbb N_{>0}$ of
	positive integers (see \eqref{eq:KLoop}) that yields decompositions of the
	transversal and orbital generating functions
	of all musical scales in
	$\varphi\big(\mathcal L (X,\mathfrak s)\big)$,
	where $\mathcal L (X,\mathfrak s)$ is the language of all admissible
	words in $X$ that start with $\mathfrak s$.
	These decompositions are
	\begin{align}
	\label{eq:decompositionT}
		\textnormal{dim}_{\mathsf T}^{\varphi (X,\mathfrak s)}(z)
		&= W^{\mathcal K(\mathfrak s)}(z) + a^{\mathcal K (\mathfrak s)}(z), \\
	\label{eq:decompositionO}
		\textnormal{dim}_{\mathsf O}^{\varphi (X,\mathfrak s)}(z)
		& = C^{\mathcal K(\mathfrak s)}(z) + b^{\mathcal K(\mathfrak s)}(z),
	\end{align}
	where the four generating functions on the right hand sides
	above are as follows:
	\begin{enumerate}
	\item \label{thm:cond1}
		$C^{\mathcal K(\mathfrak s)}(z)$ and $W^{\mathcal K(\mathfrak s)}(z)$
		are the generating functions of integer compositions
		and wheels, respectively,
		both with summands in $\mathcal K(\mathfrak s)$
		(see \eqref{eq:compOGF} and \eqref{eq:wheelsGF} evaluated at $u=1$).
	\item \label{thm:cond2}
		$a^{\mathcal K (\mathfrak s)}(z)$ is the generating function of
		aperiodic compositions, also with summands in $\mathcal K(\mathfrak s)$,
		except for the last one that belongs to
		the complement $\mathcal K(\mathfrak s)^{\textnormal{\textsf c}}$
		and is bounded above by an element of $\mathcal K(\mathfrak s)$
		(see \eqref{eq:dimTOC}).
	\item \label{thm:cond3}
		The orbital generating function $b^{\mathcal K (\mathfrak s)}(z)$ associated to the class
		represented by $a^{\mathcal K(\mathfrak s)}(z)$ is, in fact,
		the corresponding cumulative generating function
		with respect to the number of notes (see \eqref{eq:CGF}).
	\end{enumerate}
\end{theorem}

The proof follows from the definitions and the formulations
in the rest of the paper, which is organized as follows.
In section \ref{sec:ScalesCompositions}
we declare the class of all musical scales as combinatorially isomorphic
to the integer compositions, and define their modes, orbitals,
transversals, and their dimensions. 
In section \ref{sec:VertexShifts} we address shift spaces
and the classes of musical scales they define.
As first examples we mention the cases of the
Thue-Morse scales, which is a substitutive system,
and also the Fibonacci and the Fagenbaum scales, also substitutive.
We recall periodic points and zeta functions and then
focus on vertex shifts and loop systems.
Finally we settle the decompositions in the Theorem \ref{thm:main}.
Our formalism is general for all $n$-TET tuning systems,
and adapts for finite values of $n$ to
numerical procedures for exact computations,
for example when $n=12$, 
which is of most interest from 
the musical point of view.
Here, as an example to illustrate the methods in sections
\ref{sec:ScalesCompositions} and \ref{sec:VertexShifts}
that yield Theorem \ref{thm:main}, we develop
the golden mean scales.
In the last section \ref{sec:discussion} we make
final remarks and conclusions, with respect to other works,
possible generalizations, and further applications.

\section{Scales, modes, orbitals, and transversals}
\label{sec:ScalesCompositions}

\subsection{Musical scales and integer compositions}

A musical scale in 12-TET tuning system can be coded by a sequence of
integers in which each term is the number of half tones within consecutive
notes of the scale. For example, the chromatic and major scales, that are
coded with a binary alphabet $\mathcal B = \{\circ, \bullet \}$ as
$\circ\circ\circ\circ\circ\circ\circ\circ\circ\circ\circ\circ\phantom{,}\hspace{-0.18cm}$ and
$\circ\bullet\circ\bullet\circ\circ\bullet\circ\bullet\circ\bullet\circ\phantom{,}\hspace{-0.18cm}$,
correspond to $(1,1,1,1,1,1,1,1,1,1,1,1)$ and $(2,2,1,2,2,2,1)$ in the integer sequence representation, respectively.
Observe that in both cases the sums of the entries yields 12.
This is a general phenomena that, though elementary by definition,
we state formally:

\begin{theorem}
\label{thm:ScalesCompositions}
	In $n$-TET tuning system, the musical scales are in bijective
	correspondence with the set of ordered sequences of positive integers that
	add up to $n$. In other words, the set of all musical scales (in any tuning system)
	is combinatorially isomorphic to the combinatorial class of integer
	compositions.
\end{theorem}

Let $\mathbb N_{>0} \triangleq \{1,2,3,\ldots \}$ denote the set of
positive integers. Let
\begin{align}
\label{eq:comspec}
	\mathcal C \triangleq \textsc{Seq}(\mathbb N_{>0}) = \bigcup\limits_{k=0}^\infty \mathbb N_{>0}^k
\end{align}
denote the class of integer compositions,
and henceforth think of its elements as musical scales.
For every integer $n\geq 0$,
let $\mathcal C_n \subset \mathcal C$ be the compositions of $n$,
i.e. $\mathcal C_n$ denotes the set of scales in $n$-TET tuning system.
Then $C_n \triangleq \# \mathcal C_n = 2^{n-1}$, which is consistent with
the binary representation of musical scales.
Thus, the ordinary generating function (OGF) of all musical scales is the rational
function\footnote{Combinatorially, integer compositions are sequences of positive integers
(see \eqref{eq:comspec}),
and on the right hand side of \eqref{eq:compoGF} we already see the form of P\'olya's quasi-inverse operator
that corresponds to sequence constructions.}
\begin{align}
\label{eq:compoGF}
	C(z) \triangleq \sum\limits_{n=0}^\infty C_nz^n = \frac{1-z}{1-2z} .
\end{align}
For any arbitrary sequence $x=(x_1,\dots,x_k)$,
its \emph{length} is denoted by $\ell(x)\triangleq k$
(e.g., when coded by integer compositions,
the major and chromatic scales have lengths 7 and 12, respectively,
but on the other hand, both have length 12 in binary code).
More generally, for any
$\mathcal K \subseteq \mathbb N_{>0}$,
the class $\mathcal C^{\mathcal K}\subseteq \mathcal C$
of all integer compositions with summands in $\mathcal K$,
together with the length $\ell\colon\mathcal C^{\mathcal K} \to \mathbb N$
as a parameter, has bivariate generating function (BGF)
\begin{align}
\label{eq:compositionsBGF}
	C^{\mathcal K} (z,u) \triangleq \sum\limits_{n,m\geq 0} C^{\mathcal K}_{n,m} z^nu^m =
	\frac{1}{1-u\sum\limits_{k\in \mathcal K}z^k},
\end{align}
where
$C^{\mathcal K}_{n,m}
\triangleq \# \{w \in \mathcal C_n^{\mathcal K} : \ell (w) = m \}$
and $C^{\mathcal K}_n \triangleq
\#(\mathcal C^\mathcal K \cap \mathcal C_n)$. In particular,
the OGF of $\mathcal C^{\mathcal K}$ is
\begin{align}
\label{eq:compOGF}
	C^{\mathcal K} (z) \triangleq \sum\limits_{n=0}^\infty C^{\mathcal K}_n z^n =
	C^{\mathcal K} (z,1).
\end{align}

\subsection{Modes, orbitals, transversals, and their dimensions}

\begin{figure}[t]
   \centering
   \includegraphics[width=1.6in]{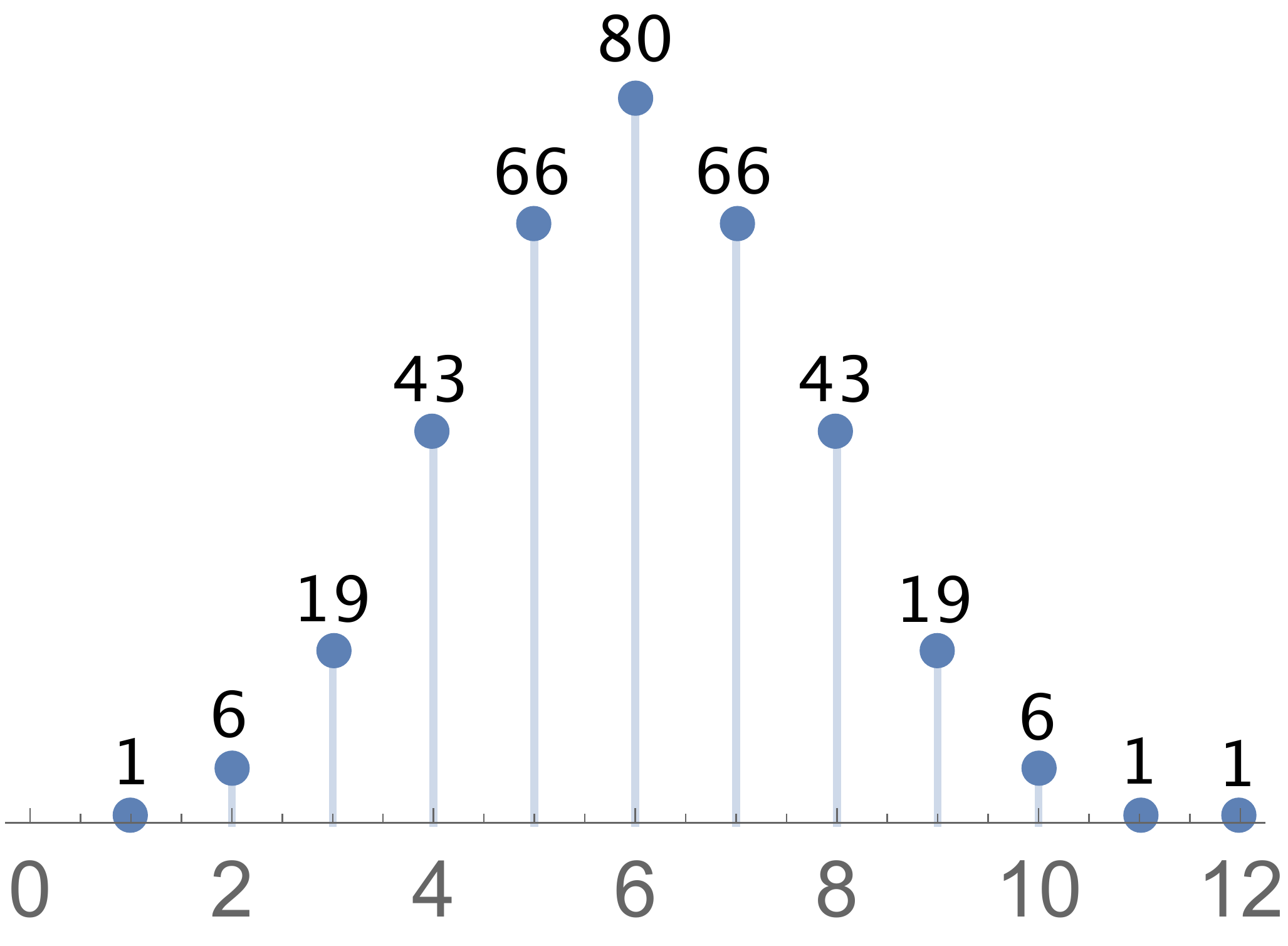} 
   \caption{Number of notes versus number of essentially different scales for the 12-TET tuning system.
   This is in general represented for every $n$ by the coefficient $[z^{n}]W(z,u)$,
   where $W(z,u)$ is the bivariate version
   \eqref{eq:wheelsGF} of $W(z)$ (when $\mathcal K = \mathbb N_{>0}$). Thus, for example,
   $[z^{12}]W(z,u) = 
   u + 6 u^2 + 19 u^3 + 43 u^4 + 66 u^5 + 80 u^6 + 66 u^7 + 43 u^8 + 
 19 u^9 + 6 u^{10} + u^{11} + u^{12}.$
   The limiting distribution of the number of notes is gaussian as $n\to \infty$.}
   \label{fig:TransversalLength}
\end{figure}

Let $\mathcal A$ be a countable alphabet
and then let
$\mathcal A^* \triangleq \bigcup_{k\geq 0} \mathcal A^k$,
where
$
	\mathcal A^k \triangleq
	\underset{k \textnormal{ times}}{\underbrace{\mathcal A \times \cdots \times \mathcal A}}.
$
Let $\alpha\colon \mathbb Z \curvearrowright \mathcal A^*$
be the cyclic left shift action induced by the combinatorial isomorphism
$\alpha \colon \mathcal A^* \to \mathcal A^*$ defined for every $w=(w_1, \ldots , w_{k})\in\mathcal A^k$
by $\alpha(w)\triangleq (w_2,\ldots, w_k, w_1) \in \mathcal A^k$, for all $k\geq 1$.
The $\alpha$-\emph{orbit} of $w\in\mathcal A^*$ is
$$
	\mathcal O_\alpha(w) \triangleq \{\alpha^j(w) \ : \ \forall \ j \in \mathbb Z \}.
$$
The set of $\alpha$-orbits forms a partition of $\mathcal A^*$ induced by
the $\alpha$-orbit equivalence relation $\overset{\alpha}{\sim}$.
The representation of musical scales
by integer compositions is such that
the $\alpha$-orbit equivalence class of an integer composition $\mathbf w\in \mathcal C$,
i.e. the elements of its $\alpha$-orbit $\mathcal O_\alpha(\mathbf w)$,
are the \emph{modes} of the corresponding scale\footnote{This is not generally the case. For example, in the binary representation of musical scales,
the $\alpha$-orbits do not always correspond to the modes of the scales.},
and thus, in this case, we write
$\textsc{modes}(\mathbf w)\triangleq \mathcal O_\alpha(\mathbf w)$.
For any subset $B\subseteq \mathcal A^*$, let
$$
	\mathcal O_\alpha (B) \triangleq \bigcup_{w\in B}\mathcal O_\alpha(w),
$$
and similarly, if $B\subseteq \mathcal C$, then we write
$\textsc{modes}(B)\triangleq \mathcal O_\alpha(B)$.
Now, since $\mathcal A^* / \overset{\alpha}{\sim}$ is the combinatorial class
of \emph{cycles} of elements of $\mathcal A$, the class of all musical scales, modulo their modes,
is the class $\mathcal W$ of \emph{cyclic} compositions of positive integers,
the so called \emph{wheels}. For example,
the diatonic wheel $(2,2,1,2,2,2,1)$ has size 12, length 7,
it is aperiodic, thus it consists of 7 modes.
Two musical scales are \emph{essentially different} if they are different as wheels.
Therefore, the OGF of all musical scales, modulo their modes, is
\begin{align}
	W(z) & \triangleq 
	\sum\limits_{k=1}^\infty W_n z^n \\
	& = \sum\limits_{k=1}^\infty \frac{\phi(k)}{k} \log \left(1-\frac{z^k}{1-z^k}\right)^{-1} \\
	& = z+2z^2+3z^3+5z^4+7z^5+13z^6+\ldots ,
\end{align}
with $\phi \colon \mathbb N_{>0} \to \mathbb N_{>0}$ the Euler totient function, that is,
$\phi(n) \triangleq \# \{k \leq n :  \gcd(n,k) = 1\}$. 
More generally, the BGF
of the class $\mathcal W^{\mathcal K}$ of wheels with summands in $\mathcal K$,
with $u$ marking the length of the wheels (i.e. the number of notes in the scales), is
\begin{align}
\label{eq:wheelsGF}
	W^{\mathcal K}(z,u) \triangleq \sum\limits_{n,m=1}^\infty W_{n,m}^{\mathcal K}z^nu^m
	= \sum\limits_{k=1}^\infty
	\frac{\phi(k)}{k} \log \frac{1}{1-\sum\limits_{j\in \mathcal K} u^kz^{jk}}.
\end{align}
For example, in 12-TET tuning system,
there are 351 essentially different musical scales, and their
distribution according
to the number of notes
is illustrated in Figure \ref{fig:TransversalLength}.

A set $A\subseteq \mathcal A^*$ is \emph{independent}
if any pair of distinct elements of $A$ belong to distinct $\alpha$-orbit
equivalence classes, i.e.
$\mathcal O_\alpha(v)\cap\mathcal O_\alpha(w)=\varnothing$ for all $v,w\in A$ with $v\neq w$.
Two subsets $A,B\subseteq \mathcal A^*$ are
\emph{mutually independent}
if $\mathcal O_\alpha(A)\cap \mathcal O_\alpha(B)=\varnothing$
(each set $A$ and $B$ may or may not be independent).
A \emph{transversal} of $A$ is a maximal independent
subset $T_A\subseteq A$.
Any nonempty set $A\neq \varnothing$
possesses at least one transversal $T_A\subseteq A$, and any two transversals of $A$ have
the same cardinality, the \emph{transversal dimension}
$$
	\textnormal{dim}_{\textsf{T}}(A) \triangleq \# T_A.
$$
Clearly, $A \subseteq \mathcal O_\alpha(A)=\mathcal O_\alpha(T)$
and $\mathcal O_\alpha(T') \subsetneq \mathcal O_\alpha(T)$ for any transversal
$T\subseteq A$ and any proper (independent) subset $T'\subsetneq T$.
Hence, if $A\subseteq \mathcal C$ is a set of integer compositions,
then the transversal dimension is the number of essentially different
scales an instrumentalist would have to learn to play any scale in $\textsc{modes}(A)$.
The \emph{orbital dimension} of $A$ is
$$
	\textnormal{dim}_{\textsf{O}}(A) \triangleq \#\mathcal O_\alpha (A).
$$
Again for subsets $A\subseteq \mathcal C$ of integer compositions,
the orbital dimension $\textnormal{dim}_{\textsf{O}}(A) = \# \textsc{modes}(A)$ is
the total number of scales that an instrumentalist can play with
the elements of (a transversal $T\subseteq A$ of) $A$
and their modes.
Thus, the orbital dimension 
$\textnormal{dim}_{\textsf O}(w) \triangleq \#\mathcal O_\alpha(w)$ of $w \in \mathcal A^*$
is bounded above by $\ell(w)$. Moreover, the former divides the later, i.e.
$\textnormal{dim}_{\textsf O}(w)|\ell(w)$, thus there is an integer $k=k(w)\geq 1$
such that $\textnormal{dim}_{\textsf O}(w) \cdot k = \ell(w)$,
and then let the \emph{period} of $w$ be defined as $\textnormal{per}(w) \triangleq k$.
If $\textnormal{per}(w)=1$, then $w$ is \emph{aperiodic}. The orbital dimension of $A$
is therefore computed with any $\alpha$-transversal $T \subseteq A$ through the equality
\begin{align}
\label{eq:orbitPeriod}
	\textnormal{dim}_{\textsf O}(A)
	= \sum\limits_{w \in T} \frac{\textnormal{dim}_{\textsf O}(w)}{\textnormal{per}(w)}.
\end{align}

\begin{example}
The set of standard Thue-Morse scales reported in
\cite{GomezNasser20} is not independent. The reader
can easily verify that the transversal and orbital dimensions
of the Thue-Morse scales is 8 and 49, respectively.
The following is a transversal of admissible (standard) Thue-Morse scales:
\footnotesize
$$
\begin{array}{|l||l|l|l|l|}
	\hline
	\ell = 6 & (3,2,1,3,1,2) \  & (3,2,1,2,3,1) \   & (3,1,3,2,1,2) \   & (3,1,2,3,1,2) \\
	\hline \hline
	 \ell = 7 & (1,3,1,2,3,1,1) \ & (1,2,3,2,1,2,1) \ & (1,3,2,1,3,1,1) \ & (2,1,2,3,2,1,1) \\
	 \hline
\end{array}
$$
\normalsize
\flushright$\Diamond$
\end{example}

\section{Symbolic dynamical scales}
\label{sec:VertexShifts}

\subsection{Shift spaces and their musical scales}
A \emph{shift space} $X\subseteq \mathcal A^\mathbb Z$
is determined by a set of \emph{forbidden blocks} $\mathcal F \subseteq \mathcal A^*$, that is, 
$X = \mathsf X_{\mathcal F}$ where
$$
	\mathsf X_{\mathcal F} \triangleq \{ x=(x_n)_{n \in \mathbb Z} \in \mathcal A^\mathbb Z:
	\forall w \in \mathcal F, \ \forall k \in \mathbb Z, \ x_{[k,k+\ell(w)-1]} \neq w \}
$$
(above, and henceforth, for any sequence $x$, let
$x_{[i,j]} \triangleq x_i\dots x_j \triangleq (x_i,\dots, x_j)$),
and is accompanied by the \emph{left shift $\mathbb Z$-action}
$\sigma  \colon \mathbb Z \curvearrowright X$ induced by the automorphism
$$
	\sigma(x)_n\triangleq x_{n+1} \ \ \forall \ x=(x_n)_{n\in \mathbb Z} \in X, \ \forall n \in \mathbb Z.
$$
The \emph{language} of a shift space $X$ is $\mathcal L(X) \triangleq \bigcup\limits_{n\geq 0} \mathcal L_n(X) \subseteq \mathcal A^*$, where
$$
	\mathcal L_n(X) \triangleq \{x_{[1,n]} \in \mathcal A^n : x\in  X \},
$$
and also, for every symbol $\mathfrak s \in \mathcal A$, let
$\mathcal L(X,\mathfrak s) \triangleq \bigcup_{n\geq 1} \mathcal L_n(X,\mathfrak s)$, where
$$
	\mathcal L_n(X,\mathfrak s) \triangleq \{x_{[1,n]} \in \mathcal L_n(X) : x_1=\mathfrak s \}.
$$
$X$ is \emph{irreducible} if for every $u,w\in\mathcal L(X)$, there exists $v \in \mathcal L(X)$
such that $uvw \in \mathcal L(X)$.

The \emph{distinguished symbol rule} $\varphi \colon \mathcal L(X) \to \mathcal C$
is defined for each $w =w_1\ldots w_{n}\in \mathcal L_n(X)$
as certain composition $\varphi (w) \in \mathcal C_n$
of $n = \ell(w)$, as follows. Let $\mathfrak s \triangleq w_1$
and then let $1=n_1<n_2<\ldots<n_{r(w)}\leq n$ be the
coordinates where $\mathfrak s$ occurs in $w$, that is, $w_{j}= \mathfrak s$
if and only if $j=n_i$ for some $i=1,\ldots , r(w)$. Then the composition of $n$
that $w$ induces has length $\ell\big(\varphi(w)\big)=r(w)$ and is defined by
\begin{align}
	\label{eq:InducedScale}
	\varphi(w) \triangleq (\underset{k_1}{\underbrace{n_2-n_1}},
	\underset{k_2}{\underbrace{n_3-n_2}}, \ldots, \underset{k_{r(w)-1}}{\underbrace{n_{r(w)}-n_{r(w)-1}}}, \underset{k_{r(w)}}{\underbrace{n- n_{r(w)}+1}}).
\end{align}

\begin{example}
	The Fibonacci and Fagenbaum binary sequences are
	defined by morphisms, like the Thue-Morse sequence (see \cite{Kurka03}):
	$$
	\begin{array}{c|c|c}
	\textnormal{Thue-Morse} & \textnormal{Fibonacci} & \textnormal{Fagenbaum} \\
	\hline
	\textnormal{\oeis{A010060}} & \textnormal{\oeis{A005614}} & \textnormal{\oeis{A035263}}\\
	\hline
	\begin{array}{l}
		\circ \mapsto \circ\bullet \\
		\bullet \mapsto \bullet \circ
	\end{array}
	&
	\begin{array}{l}
		\circ \mapsto \bullet \\
		\bullet \mapsto \bullet \circ
	\end{array}
	&
	\begin{array}{l}
		\circ \mapsto \bullet\bullet \\
		\bullet \mapsto \bullet \circ
	\end{array}
	\end{array}
	$$
	One difference though is that the Fibonacci and the Fagenbaum sequences
	are not closed under bitwise negation. The 
	transversal and orbital dimensions of the Fibonacci scales are 10 and 66,
	respectively. The following is a transversal of admissible Fibonacci scales:
	\footnotesize
	$$
	\begin{array}{|c||c||l|l|l|l|}
		\hline
		\bullet & \ell = 5 & (3,3,2,3,1) \  & (3,2,3,3,1) \   & (3,2,3,2,2) \   & (2, 3, 2, 2)  \\
		\hline \hline
		 \circ & \ell = 7 & (2,2,1,2,1,2,2) \ & (2,2,1,2,2,1,2) & & \\
		 \hline \hline
		 \circ & \ell = 8 & (1,2,1,2,2,1,2,1) \ & (1,2,2,1,2,1,2,1) \ & (1,2,2,1,2,2,1,1) \ & (2,1,2,1,2,2,1,1) \\
		 \hline
	\end{array}
	$$
	\normalsize
	Similarly, for the Fagenbaum scales the transversal and orbital dimensions
	are 6 and 28, respectively, and a transversal of admissible Fagenbaum scales is
	\footnotesize
	$$
	\begin{array}{|c||c|l||c|l||c|l|}
		\hline
		\circ & \ell = 3 & (4,4,4) \ & \ell = 4 & (2,2,4,4) \ & \ell = 5 & (2,2,4,2,2) \\
		\hline\hline
		\bullet & \ell = 7 & (2,2,2,1,1,2,2) \ & \ell = 8 & (2,2,2,1,1,2,1,1) \ & \ell = 9 & (2,1,1,2,1,1,2,1,1) \\
		 \hline
	\end{array}
	$$
	\normalsize
	These three classes of scales, the Thue-Morse, the Fibonacci, and the Fagenbaum scales,
	are pairwise independent.
	\flushright$\Diamond$
\end{example}

For every $\mathfrak s \in \mathcal A$, let
$\mathcal C^{(X, \mathfrak s)} \triangleq \varphi \big(\mathcal L(X,\mathfrak s) \big)$
and
$\mathcal C^{(X)} \triangleq \varphi\big( \mathcal L(X)\big)$, and for every $n\geq 1$,
also let
$\mathcal C_n^{(X, \mathfrak s)}  \triangleq \varphi\big(\mathcal L_n(X,\mathfrak s)\big)$
and
$\mathcal C_n^{(X)}  \triangleq \varphi\big(\mathcal L_n(X)\big)$.
Then we define the OGFs
\begin{align}
	C^{(X,\mathfrak s)}(z) \triangleq \sum\limits_{n\geq 0} C^{(X,\mathfrak s)}_nz^n
	\qquad \textnormal{ and }\qquad
	C^{(X)}(z) \triangleq \sum\limits_{n\geq 0} C^{(X)}_nz^n
\end{align}
where $C^{(X,\mathfrak s)}_n \triangleq \# \mathcal C^{(X,\mathfrak s)}_n$
and $C^{(X)}_n \triangleq \# \mathcal C^{(X)}_n$.
We are concerned with the transversal and orbital BGFs
\begin{align}
	\textnormal{dim}_{\mathsf{T}}^{(X,\mathfrak s)}(z,u)
	&\triangleq \sum\limits_{n,m\geq 0}
	\textnormal{dim}_{\mathsf{T}}\big( \mathcal C_{n,m}^{(X, \mathfrak s)} \big)z^nu^m,
	& \textnormal{dim}_{\mathsf{O}}^{(X,\mathfrak s)}(z,u)
	&\triangleq \sum\limits_{n,m\geq 0} \textnormal{dim}_{\mathsf{O}}\big(\mathcal C_{n,m}^{(X, \mathfrak s)}\big)z^nu^m, \\
	\textnormal{dim}_{\mathsf{T}}^{(X)}(z,u)
	& \triangleq \sum\limits_{n,m\geq 0} \textnormal{dim}_{\mathsf{T}}\big(\mathcal C_{n,m}^{(X)}\big)z^nu^m,
	& \textnormal{dim}_{\mathsf{O}}^{(X)}(z)
	& \triangleq \sum\limits_{n,m\geq 0} \textnormal{dim}_{\mathsf{O}}\big(\mathcal C_{n,m}^{(X)}\big)z^nu^m.
\end{align}

\subsection{Periodic points and zeta functions}
The \emph{$\sigma$-orbit} of $x \in X$ is
$\mathcal O_{\sigma}(x) \triangleq \{ \sigma^n(x) : n \in \mathbb Z\} \subseteq X$.
For every $n\geq 1$, a point $x \in X$ is \emph{$n$-periodic} if $\sigma^n(x) = x$,
and if $x$ is $n$-periodic, then there exists the \emph{minimal period} $n_x \geq 1$ of 
$x$, namely, the cardinality of its orbit $n_x \triangleq \# \mathcal O_{\sigma}(x)$, and moreover, $n_x|n$.
Let $P_n(X) \triangleq \{ x \in X : \sigma^n(x) = x \}$ and
$Q_n(x) \triangleq \{ x \in P_n(X) : n_x=n\}$ be the sets of
$n$-periodic points and minimal $n$-periodic points,
respectively, and also let $p_n(X) \triangleq \#P_n(X)$ and $q_n(X) \triangleq \# Q_n(X)$.
Recall that the relationship between $p_n(X)$ and $q_n(X)$ is through M\"obius inversion,
\begin{align}
\label{eq:Mobius}
	p_n(X) = \sum\limits_{k|n} q_k(X)
	\qquad \textnormal{ and } \qquad
	q_n(X) = \sum\limits_{k|n} \mu\left(\frac{n}{k}\right)p_k(X),
\end{align}
where $\mu \colon \mathbb N_{>0} \to \{-1,0,1\}$ is the M\"obius function defined by
$$
	\mu(n) \triangleq \left\{
	\begin{array}{ll}
	0  & \textnormal{if there exists $p\geq 2$ such that $p^2 | n$, and}\\
	(-1)^r & \textnormal{if $n=p_1 \cdots p_r$ with $p_1,p_2,\dots,p_r \geq 2$
	distinct prime numbers}.
	\end{array}
	\right.
$$
The \emph{dynamic zeta function} of $X$ is
$$
	\zeta_X(z) \triangleq \exp \left(\sum\limits_{n=1}^\infty \frac{p_n(X)}{n}z^n\right)
	= \prod\limits_{n\geq 1} \frac{1}{(1-z^n)^{q_n(X)/n}}.
$$

\subsection{Shifts of finite type, vertex shifts, and loop systems}
A \emph{shift of finite type} $X \subseteq \mathcal A^\mathbb Z$
is a shift space $X=\mathsf X_{\mathcal F}$ that can be defined by a \emph{finite}
set of forbidden blocks $\mathcal F \Subset \mathcal A^*$, and in this case 
define $m \triangleq \max \{\ell(w) : w \in \mathcal F\} - 1$ and say that
$X$ is \emph{$m$-step} (since it is always possible to find a set of $(m+1)$-blocks
$\mathcal F' \subseteq \mathcal A^{m+1}$ such that $X=\mathsf X_{\mathcal F'}$).
A \emph{vertex shift space} is a 1-step
shift of finite type. Let $X=\mathsf X_{\mathcal F} \subseteq \mathcal A^\mathbb Z$
be a vertex shift space defined by a set of forbidden 2-blocks $\mathcal F \subseteq \mathcal A^2$.
Let $A$ be the square $\{0,1\}$-matrix indexed by $\mathcal A$ and defined
by the rule $A(i,j) = 1$ if and only if $ij\notin \mathcal F$. Then $X=\widehat{\mathsf X}_A$, where
$$
	\widehat{\mathsf X}_{A} \triangleq \{x = (x_n)_{n\in \mathbb Z} \in \mathcal A ^\mathbb Z
	: \forall n \in \mathbb Z, \ A(x_n,x_{n+1}) \neq 0\}.
$$
The matrix representation of vertex shifts yields expressions that can be
useful to study transversal and orbital dimensions. For example,
the dynamic zeta function is obtained through
\begin{align}
\label{eq:zeta}
	\zeta_{\widehat{\mathsf X}_A}(z)= \frac{1}{\textnormal{det}(I-zA)}
\end{align}
(observe that $\textnormal{det}(I-zA) = z^{\#\mathcal A} \chi_A(z^{-1})$, where $\chi_A(z)$
is the characteristic polynomial of the matrix $A$).
From here we can get
$$
	p_n(\widehat{\mathsf X}_{A})
	= \left. \frac{1}{(n-1)!}\frac{d^n}{dz^n}\log \zeta_{\widehat{\mathsf X}_A}(z)\right|_{z=0}
	= \textnormal{trace}(A^n)
$$
and also $q_n(\widehat{\mathsf X}_{A})$ by M\"obius inversion \eqref{eq:Mobius}.
The following result follows.

\begin{theorem}[Transversal and orbital dimensions of languages of vertex shifts]
\label{thm:transorbitalVertex}
	The $n$th transversal dimension of the language of a vertex shift is
	\begin{align}
	\label{eq:transversalOGFSFT}
		\textnormal{dim}_{\mathsf T}\big(\mathcal L_n(\widehat{\mathsf X}_{A})\big)
		 = \sum\limits_{\substack{i,j \in \mathcal A \\ A_{j,i} = 0 }} A^{n-1}_{i,j} 
		 + \sum\limits_{k|n}\frac{q_k(\widehat{\mathsf X}_{A})}{k}
	\end{align}
	and the corresponding $n$th orbital dimension is
	\begin{align}
	\label{eq:orbitalOGFSFT}
		\textnormal{dim}_{\mathsf O}\big(\mathcal L_n(\widehat{\mathsf X}_{A})\big)
		 = n\sum\limits_{\substack{i,j \in \mathcal A \\ A_{j,i} = 0 }} A^{n-1}_{i,j} 
		 + p_n(\widehat{\mathsf X}_{A}).
	\end{align}
\end{theorem}

Now, for studying musical scales arising
from languages of vertex shift spaces
through the distinguished symbol rule, consider
the \emph{first return loop system} to a given symbol $\mathfrak s \in \mathcal A$,
defined by the OGF
\begin{align}
\label{eq:loopsystem}
	f^{(\mathfrak s)} (z) \triangleq \sum\limits_{k=1}^\infty f_k^{(\mathfrak s)}z^k,
\end{align}
where
\begin{align}
\label{eq:loopCoeff}
	f_k^{(\mathfrak s)}
	\triangleq \# \{w=w_0,\dots w_k \in \mathcal L_{k+1}(X) 
	:  w_0 w_k = \mathfrak s \neq w_j \ \ \forall \ j\neq 0,k \},
\end{align}
and that is obtained through the equation
\begin{align}
\label{eq:FirstReturn}
	1-f^{(\mathfrak s)}(z) = \frac{\zeta_{\widehat{\mathsf{X}}_B}(z)}{\zeta_{\widehat{\mathsf{X}}_{A}}(z)},
\end{align}
where $B$ is the square $\{0,1\}$-matrix indexed by $\mathcal A \setminus \{\mathfrak s\}$
and obtained from $A$ by removing the row and column
indexed by $\mathfrak s$.

\subsection{Generating functions for distinguished symbol rule on vertex shifts}
Here we proof Theorem \ref{thm:main}. Let
\begin{align}
\label{eq:KLoop}
	\mathcal K(\mathfrak s) \triangleq \{k \geq 1 : f_k^{(\mathfrak s)} \neq 0\}
\end{align}
and also denote its complement by $\mathcal K(\mathfrak s)^{\textsf c} \triangleq \mathbb N_{>0} \setminus \mathcal K(\mathfrak s)$.
According to
\eqref{eq:InducedScale} and \eqref{eq:loopCoeff}, if $w \in \mathcal L(X,\mathfrak s)$,
then $\varphi(w) = (k_1,k_2,\ldots , k_{\ell(\varphi(w))})$
is a composition of $\ell(w)$, with summands in $\mathcal K(\mathfrak s)$,
except perhaps for the last summand $k_{\ell(\varphi(w))}$.
Suppose that this is the case, that is, $k_{\ell(\varphi(w))} \in \mathcal K(\mathfrak s)^{\textsf c}$.
Since $X$ is irreducible, there exists $v \in \mathcal L(X,\mathfrak s)$
such that $v$ also ends in $\mathfrak s$ and $w$ is a prefix of $v$, that is, $v_{\ell(v)}=\mathfrak s$ and
$$
	v = wv_{[\ell(w)+1, \ell(v)]},
$$
thus $k_{\ell(\varphi(w))}$ is bounded above by an element of $\mathcal K(\mathfrak s)$.
Let
$$
	a^{\mathcal K(\mathfrak s)}(z)\triangleq \sum\limits_{n\geq 1} a_n^{\mathcal K(\mathfrak s)}z^n
$$
be the OGF of this
subclass which is described in item \ref{thm:cond2} of Theorem \ref{thm:main}.
Then
\begin{align}
\label{eq:dimTOC}
	a_{n}^{\mathcal K (\mathfrak s)}   & =
		 \sum\limits_{\substack{k \in \mathcal K(\mathfrak s)^{\textsf c} \\ \exists k'\in\mathcal K(\mathfrak s), \ k'>k}}
		C_{n-k}^{\mathcal K(\mathfrak s)} 
\end{align}
(for any subset $\mathcal K \neq \varnothing$ of positive integers,
$C_{0}^{\mathcal K} \triangleq  1$).
If we also let
$$
	b^{\mathcal K(\mathfrak s)}(z)\triangleq \sum\limits_{n\geq 1} b_n^{\mathcal K(\mathfrak s)}z^n
$$
be the OGF of the corresponding orbital,
then, by independence, there is a decomposition of
the form \eqref{eq:decompositionT} and \eqref{eq:decompositionO},
as described in Theorem \ref{thm:main},
we just need to justify that $b^{\mathcal K(\mathfrak s)}(z)$ is the cumulative generating function
of the subclass represented by $a^{\mathcal K(\mathfrak s)}(z)$, with respect to the number of notes.
This follows from the fact that the elements represented by $a^{\mathcal K(\mathfrak s)}(z)$ are
aperiodic. To be explicit, write the bivariate coefficients
\begin{align}
\label{eq:dimTransversal}
	\textnormal{dim}_{\mathsf T}\big(\mathcal C_{n,m}^{(X,\mathfrak s )}\big)  =
	W_{n,m}^{\mathcal K(\mathfrak s)} + a^{\mathcal K(\mathfrak s)}_{n,m} 
	\qquad \textnormal{ and } \qquad
	\textnormal{dim}_{\mathsf O}\big(\mathcal C_{n,m}^{(X,\mathfrak s )}\big)  =
	C_{n,m}^{\mathcal K(\mathfrak s)} + b^{\mathcal K(\mathfrak s)}_{n,m} .
\end{align}
Then, for every $n,m\geq 1$, we have
\begin{align}
\label{eq:dimTPlus}
	a_{n,m}^{\mathcal K (\mathfrak s)}   & =
		 \sum\limits_{\substack{k \in \mathcal K(\mathfrak s)^{\textsf c} \\ \exists k'\in\mathcal K(\mathfrak s), \ k'>k}} 
		C_{n-k,m-1}^{\mathcal K(\mathfrak s)} .
\end{align}
By aperiodicity, the corresponding orbital dimension is
\begin{align}
	\label{eq:dimOPlus}
	b_{n,m}^{\mathcal K (\mathfrak s)}  & =
	m\cdot a_{n,m}^{\mathcal K (\mathfrak s)}.
\end{align}
Thus, if we define
\begin{align}
\label{eq:transversalMulti}
	a^{\mathcal K (\mathfrak s)}(z,u)  \triangleq \sum\limits_{n,m\geq 1} a_{n,m}^{\mathcal K (\mathfrak s)} z^nu^m
\end{align}
and
\begin{align}
\label{eq:orbitalMulti}
	b^{\mathcal K (\mathfrak s)}(z,u)  \triangleq \sum\limits_{n,m\geq 1} b_{n,m}^{\mathcal K (\mathfrak s)} z^nu^m,
\end{align}
then we observe that
\begin{align}
\label{eq:CGF}
	b^{\mathcal K (\mathfrak s)}(z) = \frac{\partial}{\partial u}a^{\mathcal K (\mathfrak s)}(z,u)|_{u=1},
\end{align}
and in fact
\begin{align}
\label{eq:CBGF}
	b^{\mathcal K (\mathfrak s)}(z,u) = u \frac{\partial}{\partial u}a^{\mathcal K (\mathfrak s)}(z,u).
\end{align}
Hence $b^{\mathcal K (\mathfrak s)}(z)$ is the cumulative generating function of
the number of summands in the class of compositions represented by $a^{\mathcal K (\mathfrak s)}(z)$.
This settles Theorem \ref{thm:main}, and also gives
decompositions of the transversal and orbital BGFs,
$$
	\textnormal{dim}_{\mathsf T}^{\varphi(X,\mathfrak s)}(z,u)
	= W^{\mathcal K(\mathfrak s)}(z,u) + a^{\mathcal K (\mathfrak s)}(z,u)
$$
and
$$
	\textnormal{dim}_{\mathsf O}^{\varphi(X,\mathfrak s)}(z,u)
	= C^{\mathcal K(\mathfrak s)}(z,u) + b^{\mathcal K (\mathfrak s)}(z,u).
$$
To determine the transversal and orbital dimensions of the whole
set of vertex shift scales
$\mathcal C^{(X)} \triangleq \varphi\big(\mathcal L(X) \big)
= \cup_{\mathfrak s \in \mathcal A}\mathcal C^{(X,\mathfrak s)}$,
it is required to take into account the intersections between each pair of symbols,
otherwise multiple counting may occur. 
The analysis can be done one symbol at the time,
adding only new contributions to the cumulative counting.

\section{Examples}
\label{sec:golden}

\subsection{Symbolic dynamical scales of finite type}

Here we use the material from sections \ref{sec:ScalesCompositions}
and \ref{sec:VertexShifts}.

\begin{example}[Golden mean scales]
	Consider the \emph{golden mean shift}
	$X \triangleq \mathsf X_{\mathcal F}\subseteq \mathcal B^\mathbb Z$
	that is defined as the subshift that
	results from the forbidden set of blocks
	$\mathcal F = \{\bullet \bullet\}$.
	Thus $X
	 = \widehat{\mathsf{X}}_{A}$ where
	$$
		A=
		\left(
		\begin{array}{cc}
			1 & 1 \\
			1 & 0
		\end{array}
		\right).
	$$
	The matrix $A$ is the adjacency matrix of the following directed graph:
	\smallskip
	
	\centerline{\includegraphics[width=0.9in]{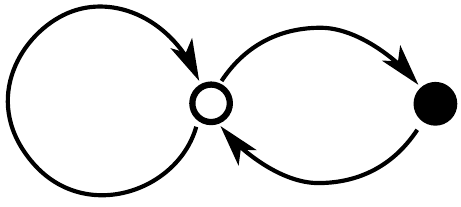}}
	
	\begin{definition}
	\label{def:fibonacci}
	For every
	$k\geq 1$, let $\{F_{n}^{(k)} \}_{n \geq 0} $ 
	be the \emph{$k$-Fibonacci sequence} defined by
	$$
		F_0^{(k)} \triangleq 1, \ \ \ F_1^{(k)} \triangleq k,
		\ \ \ \textnormal{ and } \ \ \ 
		F_{n+2}^{(k)} \triangleq F_n^{(k)}+F_{n+1}^{(k)}.
	$$
	\end{definition}

	For example, if $w_0=\circ$, $w_1=\circ\bullet$,
	and $w_{n+2}=w_{n+1}w_{n}$ for every $n\geq 0$,
	then $\lim\limits_{n\to \infty} w_n$ is the (complement of the)
	infinite Fibonacci word (\oeis{A005614}), and the sequence
	$\big( \ell(w_n) \big)_{n\geq 0} = \big(1,2,3,5,8,\ldots \big)$ is 1-Fibonacci.
	
	\smallskip

	\noindent
	\textsc{Transversals and orbitals of the golden mean language.}
	First let us illustrate the use of Theorem \ref{thm:transorbitalVertex}.
	The counting sequence of the language of the golden mean shift is
	$2$-Fibonacci, that is, $\#\mathcal L_n(X)=  F_{n}^{(2)}$ for all $n\geq 0$.
	On the other hand, the zeta function
	$$
		\zeta_X(z)=\frac{1}{1-z-z^2} = \sum\limits_{n=0}^\infty F^{(1)}_n z^n
	$$
	is the OGF of the $1$-Fibonacci sequence,
	and the periodic counting sequence
	$\big(p_n(X) = F_{n-1}^{(3)} \big)_{n\geq 1}$ is $3$-Fibonacci.
	For transversal dimensions, with reference to \eqref{eq:transversalOGFSFT},
	we see that
	$$
		\sum\limits_{\substack{i,j \in \mathcal A \\ A_{j,i} = 0 }} A^{n-1}_{i,j} =
		A^{n-1}_{\bullet,\bullet}
	$$
	represents return loops of length $n$ that begin and end at $\bullet$, but
	these are in fact sequences of \emph{first} return loops to $\bullet$.
	Using \eqref{eq:FirstReturn} and \eqref{eq:zeta} we get the OGF
	of the system of first return loops to $\bullet$, namely
	$$
		f^{(\bullet)}(z) = \frac{z^2}{1-z},
	$$
	and deduce that
	\begin{align}
		\sum\limits_{n=1}^\infty A^{n-1}_{\bullet,\bullet} z^n 
		& = z + z\frac{f^{(\bullet)}(z)}{1-f^{(\bullet)}(z)}
		= \frac{z(1-z)}{1-z-z^2} \\
		& = z + z^3 + z^4 + 2z^5 + 3 z^6 + 5 z^7 + 8 z^8
		+ 13 z^9 + \ldots , 
	\end{align}
	in particular, for $n \geq 3$, $A_{\bullet,\bullet}^{n} = F^{(2)}_{n-3}$ is 2-Fibonacci.
	Next, the sequence $\big(q_n(X)\big)_{n\geq 1}$,
	obtained from $\big(p_n(X) = F^{(3)}_{n-1}\big)_{n\geq 1}$
	by M\"obius inversion, defines a minimal periodic OGF
	\begin{align}
		q^{(X)}(z) & \triangleq \sum\limits_{n \geq 1} q_n(X) z^n \\
		& = z+z^2+z^3+z^4+2z^5+2z^6+4z^7+5z^8+8z^9 + \ldots
	\end{align}
	that corresponds to \oeis{A006206}, described as the number of
	\emph{aperiodic} binary necklaces 
	with no subsequence $\bullet \bullet$,
	excluding the necklace $\bullet$.
	The coefficients of the OGF
	\begin{align}
		\label{eq:goldenQbar}
		\overline q^{(X)}(z) & \triangleq \sum\limits_{n \geq 1} 
		\left( \sum\limits_{k | n} \frac{q_k(X)}{k} \right) z^n \\
		& = z + 2z^2+2z^3+3z^4+3z^5+5z^6+5z^7+8z^8+10z^9 + \ldots
	\end{align}
	correspond to \oeis{000358}, which is described as
	the number of necklaces
	with no subsequence $\bullet\bullet$,
	excluding the necklace $\bullet$.
	Thus the transversal OGF of the language of $X$ has the form
	$$
		\textnormal{dim}_{\textsf T}^{(X)}(z) = \frac{z(1-z)}{1-z-z^2} + \overline{q}^{(X)}(z) .
	$$
	For example, according to \eqref{eq:transversalOGFSFT},
	\begin{align}
		\textnormal{dim}_{\textsf T}\big(\mathcal L_1(X)\big) & = 1 + 1 = 2,
			& \{\circ \} \cup \{ \bullet \} \\
		\textnormal{dim}_{\textsf T}\big(\mathcal L_2(X)\big) & = 0 + 2 = 2, 
			& \{\circ\circ\phantom{.}\hspace{-0.2cm}, \circ \bullet\phantom{.}\hspace{-0.2cm} \} \\
		\textnormal{dim}_{\textsf T}\big(\mathcal L_3(X)\big) & = 1 + 2 = 3,
			&  \{\circ\circ\circ, \circ\circ\bullet \} \cup \{ \bullet \circ \bullet \} \\
		\textnormal{dim}_{\textsf T}\big(\mathcal L_4(X)\big) & = 1 + 3 = 4,
			& \{ \circ\circ\circ\circ \phantom{.}\hspace{-0.2cm}, \circ\circ\circ\bullet\phantom{.}\hspace{-0.2cm},
			\circ\bullet\circ\bullet\phantom{.}\hspace{-0.2cm}\}
			\cup \{ \bullet\circ\circ\bullet \phantom{.}\hspace{-0.2cm}\}  \\
		\textnormal{dim}_{\textsf T}\big(\mathcal L_5(X)\big) & = 2 + 3 = 5,
			& \{ \circ\circ\circ\circ\circ , \circ\circ\circ\circ\bullet,
			\circ\circ\bullet\circ\bullet\}\cup \{ \bullet\circ\circ\circ\bullet, \bullet\circ\bullet\circ\bullet \}  \\
		\textnormal{dim}_{\textsf T}\big(\mathcal L_6(X)\big) & = 3 + 5 = 8,
			& \left\{
			\begin{array}{ll}
				\circ\circ\circ\circ\circ\circ\phantom{.}\hspace{-0.2cm} , & \circ\circ\circ\circ\circ\bullet\phantom{.}\hspace{-0.2cm},\\
				\circ\circ\circ\bullet\circ\bullet\phantom{.}\hspace{-0.2cm}, & \circ\circ\bullet\circ\circ\bullet\phantom{.}\hspace{-0.2cm}, \\
				\circ\bullet\circ\bullet\circ\bullet\phantom{.}\hspace{-0.2cm}
			\end{array}
			\right\}
			\cup
			\left\{
			\begin{array}{ll}
				\bullet\circ\circ\circ\circ\bullet\phantom{.}\hspace{-0.2cm}, \\
				 \bullet\circ\circ\bullet\circ\bullet\phantom{.}\hspace{-0.2cm},  \\
				 \bullet\circ\bullet\circ\circ\bullet\phantom{.}\hspace{-0.2cm}
			\end{array}
			\right\}\\
		\vdots .
	\end{align}
	For orbital dimensions, with reference to \eqref{eq:orbitalOGFSFT}, first we see that
	\begin{align}
		\sum\limits_{n=1}^\infty  n A^{n-1}_{\bullet, \bullet} z^{n-1}
		& = \frac{d}{dz} \left(  \frac{z(1-z)}{1-z-z^2} \right) = \frac{1-2z+2z^2}{(1-z-z^2)^2} \\
		& = 1+3z^2+4z^3+10z^4+18z^5+35z^6+64z^7+117z^8+210z^9+\ldots ,
	\end{align}
	which corresponds to \oeis{A006490}.
	We already know that $p_n(X)$ is the 3-Fibonacci sequence $F^{(3)}_{n-1}$.
	Thus
	\begin{align}
		\textnormal{dim}_{\mathsf O}^{\mathcal L(X)}(z)
		& = z\frac{1-2z+2z^2}{(1-z-z^2)^2} + z\frac{1+2z}{1-z-z^2} = \frac{z(2-z-z^2-2z^3)}{(1-z-z^2)^2} \\
		& = 2z + 3z^2 + 7z^3 + 11z^4 + 21z^5 + 36z^6 + 64z^7 + 111z^8 + 193z^9 + \ldots 
	\end{align}
	(the corresponding sequence of coefficients has no record in \cite{Sloane16}).
	
	\bigskip
	
	\noindent
	\textsc{Admissible golden mean scales.}
	Now we look at the image
	$$
		\varphi\big( \mathcal L(X) \big) \subset \mathcal C
	$$
	that corresponds to the \emph{admissible golden mean scales}.
	First, using again \eqref{eq:FirstReturn} and \eqref{eq:zeta},
	we get the OGFs of the loop systems, namely
	\begin{align}
		\label{eq:goldenCirc}
		f^{(\circ)}(z) & = z+z^2 \\
		\label{eq:goldenBullet}
		f^{(\bullet)}(z) & = \frac{z^2}{1-z}=z^2+z^3+z^4+\ldots .
	\end{align}
	Then the $\circ$-admissible golden mean scales $\mathcal C^{(X,\circ)}\triangleq \varphi \big(\mathcal L(X,\circ) \big)$
	are all the integer compositions
	with summands in $\mathcal K(\circ) = \{1,2\}$ (see \eqref{eq:goldenCirc}), that is, all the scales with no more than two tone
	measures of difference between consecutive notes, as expected
	(in this case we have $a^{\mathcal K(\circ)}(z,u) = 0$,
	and thus also $b^{\mathcal K(\circ)}(z,u) = 0$,
	because any element in $\mathcal K(\circ)^{\textsf c}$ is not bounded above by any
	element of $\mathcal K(\circ)$).
	The corresponding OGF for this class of integer
	compositions, according to \eqref{eq:compOGF}, is
	\begin{align}
	\label{eq:Ccirc}
		C^{(X,\circ)}(z) = C^{\mathcal K(\circ)}(z)
		= \frac{1}{1-z-z^2},
	\end{align}
	thus $C^{(X,\circ)}_{n} = F_{n}^{(1)}$ is 1-Fibonacci.
	For example, for 12-TET tuning system,
	the total number of $\circ$-admissible golden mean scales is
	$$
		C^{(X,\circ)}_{12}=F_{12}^{(1)} = 233.
	$$
	The set $\mathcal C^{(X,\circ)}_{12}$ is too large to list.
	With the purpose of having a smaller context that helps illustrating and verifying
	the claims, let us imagine, for instance, that we are doing scales
	over a small set of notes, say over 5 notes (e.g. over a pentatonic scale).
	Combinatorially, the model is that of a 5-TET tuning system. We thus have
	\begin{align}
	\label{eq:gmcircle}
		C^{(X,\circ)}_{5}=F_{5}^{(1)} = 8,
		\qquad
		\mathcal C_5^{(X,\circ)} = 
		\left\{
		\begin{array}{llll}
			\overset{(1,1,1,1,1)}{\circ\circ\circ\circ\circ},
			\overset{(1,1,1,2)}{\circ\circ\circ\circ\bullet},
			\overset{(1,1,2,1)}{\circ\circ\circ\bullet\circ},
			\overset{(1,2,1,1)}{\circ\circ\bullet\circ\circ}, \\ \\
			\overset{(2,1,1,1)}{\circ\bullet\circ\circ\circ},
			\overset{(1,2,2)}{\circ\circ\bullet\circ\bullet},
			\overset{(2,1,2)}{\circ\bullet\circ\circ\bullet},
			\overset{(2,2,1)}{\circ\bullet\circ\bullet\circ}
		\end{array}
		\right\}.
	\end{align}	
	We also have the bivariate version of \eqref{eq:Ccirc},
	with $u$ marking the number of notes, namely
	$$
		C^{(X,\circ)}(z,u) = C^{\mathcal K(\circ)}(z,u)
		= \frac{1}{1-uz-uz^2} .
	$$

	For the $\bullet$-admissible
	golden mean scales
	$\mathcal C^{(X,\bullet)}\triangleq \varphi \big(\mathcal L(X,\bullet) \big)$,
	first observe that the class of integer compositions with summands
	in $\mathcal K(\bullet) = \{2,3,4,\ldots\}$ (see \eqref{eq:goldenBullet}) has
	OGF
	\begin{align}
	\label{eq:bulletgf}
		C^{\mathcal K(\bullet)}(z) = \frac{1-z}{1-z-z^2}.
	\end{align}
	The bivariate version of $C^{\mathcal K(\bullet)}(z)$,
	with the variable $u$ marking the number of notes, is
	$$
		C^{\mathcal K(\bullet)}(z,u) = \frac{1-z}{1-z-uz^2}.
	$$
	In addition, in this case, the last summand in the elements of $\mathcal C^{(X,\bullet)}$
	\emph{is} allowed to be $1 \notin \mathcal K(\bullet)$ 
	(for example, in 12-TET tuning system, the binary admisible 12-block
	$\bullet\circ\bullet\circ\bullet\circ\bullet\circ\bullet\circ\circ\bullet\phantom {.}\hspace{-0.2cm}$
	yields the $\bullet$-admissible golden mean scale $(2,2,2,2,3,1)$).
	Therefore, the corresponding generating function is
	\begin{align}
	\label{eq:Cbullet}
		C^{(X,\bullet)}(z) = C^{\mathcal K(\bullet)}(z) + zC^{\mathcal K(\bullet)}(z)
		= \frac{1-z^2}{1-z- z^2},
	\end{align}
	and we also get its bivariate version,
	$$
		C^{(X,\bullet)}(z,u) = C^{\mathcal K(\bullet)}(z,u) + uzC^{\mathcal K(\bullet)}(z,u)
		= \frac{(1+uz)(1-z)}{1-z- uz^2}.
	$$
	Thus $C^{(X,\bullet)}_{0}=1$ 
	and for $n\geq 1$ the sequence of coefficients
	$C^{(X,\bullet)}_{n} = F_{n-1}^{(1)}$ is 1-Fibonacci.
	For example, for the 12-TET and 5-TET tuning systems we have
	\begin{align}
		C^{(X,\bullet)}_{12} & =F^{(1)}_{11}= 144, \\
	\label{eq:gmbullet}
		C^{(X,\bullet)}_{5} & =F^{(1)}_{4}= 5,
		\qquad
		\mathcal C^{(X,\bullet)}_{5}  =
		\left\{
			\overset{(5)}{\bullet\circ\circ\circ\circ},
			\overset{(4,1)}{\bullet\circ\circ\circ\bullet},
			\overset{(3,2)}{\bullet\circ\circ\bullet\circ},
			\overset{(2,3)}{\bullet\circ\bullet\circ\circ},
			\overset{(2,2,1)}{\bullet\circ\bullet\circ\bullet} 
		\right\}.
	\end{align}
	
	\begin{remark}
	\label{rem:inter}
		The only elements that are both $\circ$-admissible and $\bullet$-admissible
		golden mean scales are the compositions of even $n=2d$ and odd $n=2d+1$ integers
		of the form
		\begin{align}
		\label{eq:intersection}
			(\underset{d\textnormal{-times}}{\underbrace{2,2,\dots,2}})
			\qquad \textnormal{ and } \qquad
			(\underset{d\textnormal{-times}}{\underbrace{2,2,\dots,2}},1).
		\end{align}
	\end{remark}
	Thus, combining \eqref{eq:Ccirc} and  \eqref{eq:Cbullet}, we conclude that
	the OGF of \emph{admissible} golden mean scales is
	\begin{align}
		C^{(X)}(z) & = C^{(X,\circ)}(z) + C^{(X, \bullet)}(z) - \frac{1+z}{1-z^2}
		= \frac{1-z+z^3}{1-2z+z^3} \\
		& = 1 + z + 2 z^2 + 4 z^3 + 7 z^4 + 12 z^5 + 20 z^6 + 33 z^7
		+ 54 z^8 + 88 z^9 + 143 z^{10}+\ldots 
	\end{align}
	(the corresponding coefficients are, essentially, \oeis{A000071}), and we
	also get its bivariate version,
	$$
		C^{(X)}(z,u) 
		= C^{(X,\circ)}(z,u) + C^{(X, \bullet)}(z,u) - \frac{1+uz}{1-uz^2}.
	$$
	Then, for every $n\geq 1$ we have $C^{(X)}_n = F_{n}^{(1)}+F_{n-1}^{(1)}-1$.
	For example, the number of admissible golden means scales
	in 12-TET and 5-TET tuning systems are (for the later see
	\eqref{eq:gmcircle} and \eqref{eq:gmbullet})
	\begin{align}
		C^{(X)}_{12} &= 233 + 144 - 1 = 376, \\
		C^{(X)}_{5} & = 8 + 5 - 1 = 12,
		\qquad \qquad
		\mathcal C_5^{(X)} = 
		\left\{
		\begin{array}{llll}
			\overset{(1,1,1,1,1)}{\circ\circ\circ\circ\circ},
			\overset{(1,1,1,2)}{\circ\circ\circ\circ\bullet},
			\overset{(1,1,2,1)}{\circ\circ\circ\bullet\circ},
			\overset{(1,2,1,1)}{\circ\circ\bullet\circ\circ}, \\ \\
			\overset{(2,1,1,1)}{\circ\bullet\circ\circ\circ},
			\overset{(1,2,2)}{\circ\circ\bullet\circ\bullet},
			\overset{(2,1,2)}{\circ\bullet\circ\circ\bullet},
			\overset{(2,2,1)}{\circ\bullet\circ\bullet\circ}, \\ \\
			\overset{(5)}{\bullet\circ\circ\circ\circ},
			\overset{(4,1)}{\bullet\circ\circ\circ\bullet},
			\overset{(3,2)}{\bullet\circ\circ\bullet\circ},
			\overset{(2,3)}{\bullet\circ\bullet\circ\circ}
		\end{array}
		\right\}.
	\end{align}
	
	\noindent
	\textsc{Transversal and orbital dimensions of golden mean scales.}
	First look at the case when $\mathfrak s = \circ$.
	With the ordinary form $W^{\mathcal K}(z) \triangleq W^{\mathcal K}(z,1)$ of \eqref{eq:wheelsGF},
	we obtain the first summand in the right hand side of \eqref{eq:decompositionT},
	$$
		W^{\mathcal K(\circ)}(z) = 
		z + 2 z^2 + 2 z^3 + 3 z^4 + 3 z^5 + 5 z^6 + 5 z^7 + 8 z^8 + 10 z^9
		+15 z^{10} + 19 z^{11} + 31 z^{12} + \ldots
	$$
	with coefficients forming again the sequence \oeis{A000358} \
	(that is, in this case, we have $W^{\mathcal K(\circ)}(z) = \overline q^{(X)}(z)$,
	see \eqref{eq:goldenQbar}). Since $a^{\mathcal K(\circ)}(z)=b^{\mathcal K(\circ)}(z)=0$,
	\begin{align}
	\label{eq:transcircgf}
		\textnormal{dim}_{\textsf T}^{\varphi(X,\circ)}(z) = W^{\mathcal K(\circ)}(z)
	\end{align}
	and
	\begin{align}
	\label{eq:orbitcircgf}
		\textnormal{dim}_{\textsf O}^{\varphi(X,\circ)}(z) = C^{\mathcal K(\circ)}(z)
	\end{align}
	(see \eqref{eq:Ccirc}). In particular, for the 12-TET and 5-TET tuning system,
	the transversal dimensions are
	\begin{align}
		\textnormal{dim}_{\textsf T}\big(\varphi\big(\mathcal L_{12}(X,\circ) \big)\big)
		& = W_{12}^{\mathcal K(\circ)}(z) = 31 \\
	\label{eq:transdim5circ}
		\textnormal{dim}_{\textsf T}\big(\varphi\big(\mathcal L_{5}(X,\circ) \big)\big)
		& = W_{5}^{\mathcal K(\circ)}(z) = 3,
		\qquad 
		T_5^{(X,\circ)} \triangleq
		\left\{
		\begin{array}{llll}
			\overset{(1,1,1,1,1)}{\circ\circ\circ\circ\circ},
			\overset{(1,1,1,2)}{\circ\circ\circ\circ\bullet},
			\overset{(1,2,2)}{\circ\circ\bullet\circ\bullet}
		\end{array}
		\right\}
	\end{align}
	and the corresponding orbital dimensions are
	\begin{align}
		\textnormal{dim}_{\textsf O}\big(\varphi\big(\mathcal L_{12}(X,\circ) \big)\big)
		& = C_{12}^{\mathcal K(\circ)}(z) = 233 \\
		\textnormal{dim}_{\textsf O}\big(\varphi\big(\mathcal L_{5}(X,\circ) \big)\big)
		& = C_{5}^{\mathcal K(\circ)}(z) = 8
		\qquad \textnormal{(see again \eqref{eq:gmcircle}).}
	\end{align}
	
	Now suppose that $\mathfrak s = \bullet$ and proceed similarly.
	For the first summand in the right hand side of \eqref{eq:dimTransversal},
	use \eqref{eq:wheelsGF} to obtain the OGF
	$$
		W^{\mathcal K(\bullet)}(z) = 
		z^2 + z^3 + 2 z^4 + 2 z^5 + 4 z^6 + 4 z^7 + 7 z^8 + 9 z^9
		+ 14 z^{10} + 18 z^{11} + 30 z^{12}+ \ldots
	$$
	with coefficients essentially forming the sequence \oeis{A032190},
	which is already described as the number of cyclic compositions of $n$ into parts $\geq 2$.
	For example, for the 12-TET and 5-TET tuning system, we have
	\begin{align}
		W_{12}^{\mathcal K(\bullet)} & = 30, \\
		W_{5}^{\mathcal K(\bullet)} & =2 ,
		\qquad
		\mathcal W_5^{\mathcal K(\bullet)}
		=\{\overset{(5)}{\bullet\circ\circ\circ\circ}, \overset{(3,2)}{\bullet\circ\circ\bullet\circ}\}.
	\end{align}
	Next, from \eqref{eq:dimTOC} we get
	$$
		a_n^{\mathcal K(\bullet)} \triangleq \sum\limits_{\substack{k \notin \mathcal K(\bullet)\\ \exists k' \in \mathcal K(\bullet), \ k' > k}}
	       C_{n-k}^{\mathcal K(\bullet)}
	       =  
	       C_{n-1}^{\mathcal K(\bullet)} 
	$$
	and thus $a_1^{\mathcal K(\bullet)}=1$, $a_2^{\mathcal K(\bullet)}=0$ and 
	$a_{n+3}^{\mathcal K(\bullet)}=F_n^{(1)}$ is the 1-Fibonacci sequence. Hence
	\begin{align}
	\label{eq:abulletgf}
		a^{\mathcal K(\bullet)}(z) & = \frac{z-z^2}{1-z-z^2} \\
		& = z + z^3 + z^4 + 2 z^5 + 3 z^6 + 5 z^7 + 8 z^8
		+ 13 z^9 + 21 z^{10} + 34 z^{11} + 55 z^{12} + \ldots .
	\end{align}
	For instance, in 12-TET and 5-TET tuning systems, we have
	\begin{align}
		a_{12}^{\mathcal K (\bullet)} & = 55 \\
	\label{eq:a5}
		a_{5}^{\mathcal K (\bullet)} & = 2,
		\qquad
		\{
		\overset{(4,1)}{\bullet\circ\circ\circ\bullet},
		\overset{(2,2,1)}{\bullet\circ\bullet\circ\bullet}
		\} .
	\end{align}
	Thus \eqref{eq:decompositionT} in Theorem \ref{thm:main} yields
	the transversal OGF
	\begin{align}
	\label{eq:transbullet}
		\textnormal{dim}_{\textsf T}^{\varphi (X,\bullet)}(z) = W^{\mathcal K(\bullet)}(z) +
		\frac{z-z^2}{1-z-z^2}.
	\end{align}
	For example, in 12-TET and 5-TET tuning system, we have
	\begin{align}
		\textnormal{dim}_{\textsf T}\big(\varphi \big(\mathcal L_{12}(X,\bullet)\big)\big)
		& = 30+55 = 85,\\
	\label{eq:transdim5bullet}
		\textnormal{dim}_{\textsf T}\big(\varphi \big(\mathcal L_{5}(X,\bullet)\big)\big)
		& = 2+2 = 4,
		\qquad
		T_5^{(X,\bullet)} = \{
		\overset{(5)}{\bullet\circ\circ\circ\circ},
		\overset{(3,2)}{\bullet\circ\circ\bullet\circ},
		\overset{(4,1)}{\bullet\circ\circ\circ\bullet},
		\overset{(2,2,1)}{\bullet\circ\bullet\circ\bullet}
		\} .
	\end{align}
	
	For orbital dimensions, first we have
	\begin{align}
	\label{eq:CbulletGF}
		C^{\mathcal K(\bullet)}(z) & = \frac{1-z}{1-z-z^2} \\
		& = 1 + z^2 + z^3 + 2 z^4 + 3 z^5 + 5 z^6 + 8 z^7 + 13 z^8 + 21 z^9
		+ 34 z^{10} + \ldots . 
	\end{align}
	For example, in 12-TET and 5-TET tuning systems, we get
	\begin{align}
	\label{eq:Cbullet12}
		C_{12}^{\mathcal K(\bullet)} & = 89, \\
	\label{eq:Cbullet5}		
		C_{5}^{\mathcal K(\bullet)} & = 3,
		\qquad
		\mathcal C_{5}^{\mathcal K(\bullet)} = 
		\{
		\overset{(5)}{\bullet\circ\circ\circ\circ},
		\overset{(3,2)}{\bullet\circ\circ\bullet\circ},
		\overset{(2,3)}{\bullet\circ\bullet\circ\circ}
		\}.
	\end{align}
	Next, according to item \ref{thm:cond3} in Theorem \ref{thm:main},
	we need the cumulative generating function
	of the class represented by $a^{\mathcal K(\bullet)}(z)$,
	with respect to the number of notes.
	From \eqref{eq:transversalMulti},
	$$
		a_{n,m}^{\mathcal K(\bullet)}
		= \sum\limits_{\substack{k \notin \mathcal K(\bullet)\\ \exists k' \in \mathcal K(\bullet), \ k' > k}}
	       C_{n-k,m}^{\mathcal K(\bullet)}
	       =  
	       C_{n-1,m}^{\mathcal K(\bullet)},
	$$
	thus, using \eqref{eq:compositionsBGF}, we get
	$$
		a^{\mathcal K(\bullet)}(z,u) = uz C^{\mathcal K(\bullet)}(z,u) = \frac{uz-uz^2}{1-z-uz^2}.
	$$
	Hence, from \eqref{eq:CBGF} we get
	\begin{align}
		b^{\mathcal K(\bullet)}(z,u) & = u\frac{\partial}{\partial u}a^{\mathcal K(\bullet)}(z,u) =
		\frac{uz(1-z)^2}{(1-z-uz^2)^2} ,
	\end{align}
	and from \eqref{eq:CGF} we obtain
	\begin{align}
	\label{eq:borbitbib}
		b^{\mathcal K(\bullet)}(z) & = \frac{z(1-z)^2}{(1-z-z^2)^2} \\
		& = z + 2 z^3 + 2 z^4 + 5 z^5 + 8 z^6 + 15 z^7 + 26 z^8 + 46 z^9
		+ 80 z^{10} + \ldots 
	\end{align}
	(the coefficients are \oeis{A006367}).
	For example, in the 12-TET and 5-TET tuning system (for the later see \eqref{eq:a5}),
	we get
	\begin{align}
	\label{eq:bbullet12}
		b_{12}^{\mathcal K(\bullet)} & = 240, \\
	\label{eq:bbullet5}
		b_5^{\mathcal K(\bullet)} & = 5,
		\qquad
		\overset{\mathcal O_\alpha ((4,1))}{\overbrace{\{ (4,1),(1,4) \}}} \cup
		\overset{\mathcal O_\alpha ((2,2,1))}{\overbrace{\{ (2,2,1), (2,1,2), (1,2,2)\}}} .
	\end{align}
	Hence, \eqref{eq:decompositionO} in Theorem \ref{thm:main},
	together with \eqref{eq:CbulletGF} and \eqref{eq:borbitbib}, yield
	\begin{align}
	\label{eq:orbitbullet}
		\textnormal{dim}_{\textsf O}^{\varphi(X,\bullet)}(z) & = 
		C^{\mathcal K(\bullet)}(z) + b^{\mathcal K(\bullet)}(z)
		= \frac{(1-z)(1-2z^2)}{(1-z-z^2)^2} \\
		& = 1 + z + z^2 + 3 z^3 + 4 z^4 + 8 z^5 + 13 z^6 + 23 z^7 + 39 z^8 + 67 z^9 + \ldots,
	\end{align}
	which corresponds to \oeis{A206268}, described as
	the number of compositions with at most one 1.
	For example, in the 12-TET and 5-TET tuning system
	(see \eqref{eq:Cbullet12}, \eqref{eq:Cbullet5}, \eqref{eq:bbullet12}, and \eqref{eq:bbullet5}),
	\begin{align}
		\textnormal{dim}_{\textsf O}\big(\varphi \big( \mathcal L_{12}(X,\bullet \big) \big)
		& = 89 + 240 = 329, \\
		\textnormal{dim}_{\textsf O}\big(\varphi \big( \mathcal L_{5}(X,\bullet \big) \big)
		& = 3 + 5 = 8,
		\qquad 
		\mathcal O_\alpha \big(\varphi \big( \mathcal L_{5}(X,\bullet \big) \big) =
		\left\{
		\begin{array}{ccc}
			(5), & (3,2), & (2,3) ,\\ (4,1), & (1,4) ,\\  ( 2,2,1), & (2,1,2), & (1,2,2)
		\end{array}
		\right\}.
	\end{align}
	
	For a global transversal, according to remark \ref{rem:inter},
	with \eqref{eq:transcircgf} and \eqref{eq:transbullet} we get 
	\begin{align}
		\textnormal{dim}_{\textsf T}^{\varphi(X)}(z) & = 
		\textnormal{dim}_{\textsf T}^{\varphi(X,\circ)}(z) +
		\textnormal{dim}_{\textsf T}^{\varphi(X,\bullet)}(z) -
		\textnormal{dim}_{\textsf T}^{\varphi(X,\circ) \cap \varphi(X,\bullet)}(z) \\
		& = W^{\mathcal K(\circ)}(z) + W^{\mathcal K(\bullet)}(z) + \frac{z-z^2}{1-z-z^2} - \frac{z}{1-z} \\
		&= z + 2 z^2 + 3 z^3 + 5 z^4 + 6 z^5 + 11 z^6 + 13 z^7 + 22 z^8
		+ 31 z^9 +\ldots 
	\end{align}
	(the corresponding sequence of coefficients has no record in \cite{Sloane16}).
	For example, for the 12-TET and 5-TET tuning system, the transversal dimensions of
	the golden mean scales are (for the later see (\eqref{eq:transdim5circ}) and \eqref{eq:transdim5bullet})
	\begin{align}
		\textnormal{dim}_{\textsf T}\big(\varphi\big(\mathcal L_{12} (X)\big) \big) & = 115 \\
		\textnormal{dim}_{\textsf T}\big(\varphi\big(\mathcal L_{5} (X)\big) \big) & = 6,
		\qquad
		T_5^{(X,\circ)} \cup T_5^{(X,\bullet)} \triangleq
		\left\{
		\begin{array}{llll}
			\overset{(1,1,1,1,1)}{\circ\circ\circ\circ\circ},
			\overset{(1,1,1,2)}{\circ\circ\circ\circ\bullet},
			\overset{(1,2,2)}{\circ\circ\bullet\circ\bullet}, \\
			\overset{(5)}{\bullet\circ\circ\circ\circ},
			\overset{(3,2)}{\bullet\circ\circ\bullet\circ},
			\overset{(4,1)}{\bullet\circ\circ\circ\bullet}
		\end{array}
		\right\}
	\end{align}
	($(1,2,2)$ and $(2,2,1)$
	are equal as wheels).
	
	Finally, for the global orbital,
	now we use \eqref{eq:orbitcircgf} and \eqref{eq:orbitbullet},
	but first observe that, according to remark \ref{rem:inter},
	in addition to the empty composition, there are two kinds
	of compositions in the intersection (see \eqref{eq:intersection}):
	one kind is formed by
	compositions of even integers $n=2d$ of length $d$ and period 1,
	and the other kind is formed by compositions of odd $n=2d+1$ of length
	$d + 1$ that are aperiodic. We conclude that
	\begin{align}
		\textnormal{dim}_{\textsf O}^{\varphi(X)}(z) & = 
		\textnormal{dim}_{\textsf O}^{\varphi(X,\circ)}(z) +
		\textnormal{dim}_{\textsf O}^{\varphi(X,\bullet)}(z) -
		\textnormal{dim}_{\textsf O}^{\varphi(X,\circ) \cap \varphi(X,\bullet)}(z) -1\\
		& = \frac{1}{1-z-z^2} +  \frac{(1-z)(1-2z^2)}{(1-z-z^2)^2}
		- \frac{z^2}{1-z^2} - \sum\limits_{d=0}^\infty (d+1)z^{2d+1} -1 \\
		& = \frac{1-z-3z^2+3 z^3 + 4 z^4 - 5 z^5 - 2 z^6 + 2 z^7}{(1-z^2)^2(1-z-z^2)^2} \\
		& = 1 + z + 2 z^2 + 4 z^3 + 8 z^4 + 13 z^5 + 25 z^6 + 40 z^7 + 72 z^8 + 117 z^9 +\ldots .
	\end{align}
	For example, in the 12-TET and 5-TET tuning system, we get
	\begin{align}
		\textnormal{dim}_{\textsf O}\big(\varphi\big(\mathcal L_{12} (X)\big) \big) & = 561 \\
		\textnormal{dim}_{\textsf O}\big(\varphi\big(\mathcal L_{5} (X)\big) \big) & = 13,
		\qquad
		\textsc{modes}\big(\mathcal L_5(X)\big) =
		\left\{
		\begin{array}{llll}
			(1,1,1,1,1), (1,1,1,2), (1,1,2,1), \\
			(1,2,1,1),(2,1,1,1), (1,2,2), \\
			(2,2,1),(2,1,2), (5), (3,2), \\
			(2,3),(4,1),(1,4)
		\end{array}
		\right\}.
	\end{align}
	\flushright{$\Diamond$}
\end{example}

\section{Related works and discussion}
\label{sec:discussion}

We have seen a general method to deduce transversal and orbital generating functions
of classes of musical scales induced by vertex shift spaces through the
distinguished symbol rule. The more general problem for SFTs yields the
``distinguished set of symbols rule'' for vertex shifts, and instead of a single
power series like \eqref{eq:loopsystem}, we get a \emph{matrix} of formal power series
(further details will be treated in a separate work).
We have found approaches that use generating
functions at least in the texts \cite{Jedrzejewski06, Benson07}, and learned that
enumeration problems in music go back at least to the works of
\cite{Reiner85, Read97,Fripertinger99}.
Our methods can serve to complement and generalize several other works that
address characterizations and classification of musical scales like
\cite{Nuno20,Kozyra14}, also
octave subdivisions \cite{HearneMilneDean19}, 
optimal spelling of pitches of musical scales \cite{BoraTezelVahaplar19},
tuning systems other that 12-TET \cite{HearneMilneDean19},
scales and constraint programming \cite{Hooker16},
modular arithmetic sequences and scales \cite{Amiot15},
algebras of periodic rhythms and scales \cite{AmiotSethares11},
formalisms to generate pure-tone systems that best approximate 
modulation/transposition properties of equal-tempered scales
\cite{KrantzDouthett11},
tuning systems and modes \cite{GarmendiaNavarro95}, etc.
Moreover, other combinatorial classes,
such as non-crossing configurations 
\cite{FlajoletNoy99} like dissections of
polygons and RNA secondary structures
\cite{Gomez15}, can be incorporated
to complement works that address constructions of
musical scales like \cite{MilneBulgerHerff15}.
There are many references that address the theory of musical scales
that are relevant to our work, like the fundamentals \cite{Forte73,Slonimsky47},
from the point of view of mathematics inclusive  \cite{Isola16, LindleyTurner93},
several of which are related to combinatorics on words
\cite{AlloucheJohnson18, ClampittNoll18, Brlek18, AbdallahGoldMarsden16}.

In our arguments, a key ingredient has been the use of first return loop systems,
which arise in the study of classification problems
of Markov shifts \cite{Gomez03,BoyleBuzziGomez06}.
In fact, studying music theory in contexts of dynamical systems has
been an active area of research, for example \cite{Amiot20}, see also \cite{TymoczkoYust19}.
Furthermore, the results presented here can serve as a basis to
adapt other related areas of mathematics in music, such as thermodynamic formalism
and random environments \cite{BarbieriGomezMarcusTaati20}
(e.g. to compute (relative) partition functions).



\section*{ORCID iDs}

\noindent
Ricardo G\'omez A\'iza: \color{blue}\href{https://orcid.org/0000-0002-2614-8519}{https://orcid.org/0000-0002-2614-8519}\color{black}



\section*{Acknowledgements}
\addcontentsline{toc}{section}{Acknowledgements}
I thank Doug Lind for pointing out to us the reference \cite{Benson07}.

\section*{Funding}
\addcontentsline{toc}{section}{Funding}

This work was supported by DGAPA-PAPIIT project IN107718.

\bibliographystyle{plain}
\bibliography{GomezSDSBibTexDatabase}

\addcontentsline{toc}{section}{References}

\end{document}